\documentclass[11pt, a4paper]{article}
\usepackage{amsfonts}
\usepackage{mathrsfs}
\usepackage{amssymb}

\usepackage{latexsym}
\usepackage{graphicx}
\newtheorem{theorem}{Theorem}

\newenvironment{proof}{\noindent {\bf Proof.}}

\long\def\longdelete#1{} \baselineskip 21 pt  

\title{{\bf On the number of SDRs of a valued  {\boldmath $(t,n)$}-family}\thanks{Supported in part by National Natural
Science Foundation of China (No. 10871166) and Shanghai Leading
Academic Discipline Project (No. B407).}}
\author{Dawei He\thanks{E-mail: davidecnu@gmail.com}\ \ and \ \ Changhong Lu\thanks{E-mail: chlu@math.ecnu.edu.cn}  \\
   Department of Mathematics,\\
                East China Normal University,\\
                Shanghai 200241, P. R. China}

\date{June 30, 2010} 
\begin{document}
\maketitle

\begin{abstract}
A system of distinct representatives (SDR) of a family $F = (A_1,
\cdots, A_n)$ is a sequence $(x_1, \cdots, x_n)$ of $n$ distinct
elements with $x_i \in A_i$ for $1 \le i \le n$. Let $N(F)$ denote
the number of SDRs of a family $F$; two SDRs are considered distinct
if they are different in at least one component. For a nonnegative
integer $t$, a family $F=(A_1,\cdots,A_n)$ is called a
$(t,n)$-family if the union of any $k\ge 1$ sets in the family
contains at least $k+t$ elements. The famous Hall's Theorem says
that $N(F)\ge 1$ if and only if $F$ is a $(0,n)$-family. Denote by
$M(t,n)$ the minimum number of SDRs in a $(t,n)$-family. The problem
of determining $M(t,n)$ and those families containing exactly
$M(t,n)$ SDRs was first raised by Chang [European J. Combin.{\bf
10}(1989), 231-234]. He solved the cases when $0\le t\le 2$ and gave
a conjecture for $t\ge 3$. In this paper, we solve the conjecture.
In fact, we get a more general result for  so-called valued
$(t,n)$-family.

\end{abstract}




\section{Introduction}

A system of distinct representatives (SDR) of a family $F = (A_1,
\cdots, A_n)$ is a sequence $(x_1, \cdots, x_n)$ of $n$ distinct
elements with $x_i \in A_i$ for $1 \le i \le n$.  The famous Hall's
theorem \cite{phall1935} tell us that  a family has a SDR if and
only if the union of any $k\ge 1$ sets of this family contains at
least $k$ elements. Several quantative refinements of the Hall's
theorem were given in \cite{mhall1948, mirsky1971, rado1967}. Their
results are all under the assumption of Hall's condition plus some
extra conditions on the cardinalities of $A_i$'s.

Chang \cite{chang1989} extends Hall's theorem as follows: let $t$ be
a nonnegative integer.  A family $F=(A_1,\cdots,A_n)$ is called a
{\it $(t,n)$-family} if $|\bigcup\limits_{i \in I}A_i|\ge |I| + t$
holds for any non-empty subset $I \subseteq \{1, \cdots, n\}$.
Denote by $N(F)$  the number of SDRs of a family $F$. Let $M(t,n) =
\min \{N(F)~|~F$ is a $(t,n)$-family$\}$. Hall's theorem says that
$M(0,n)\ge 1$.  In fact,  it is easy to know that $M(0,n)=1$.  Chang
\cite{chang1989} proved that $M(1,n)=n+1$ and $M(2,n)=n^2+n+1$. He
also determined all $(t,n)$-families $F$ with $N(F)=M(t,n)$ for
$t=0,1,2$. Consider the $(t,n)$-family $F^{*} = (A_{1}^{*}, \cdots,
A_{n}^{*})$, where $A_{i}^{*} = \{i, n+1, \cdots, n+t\}$ for $1 \le
i \le n$. Then,

$$ N(F^{*})=U(t,n)=\sum\limits_{j = 0}^{t} {t \choose j} {n \choose j}j!.$$
Chang\cite{chang1989} has shown  that $F^*$ as above is the only
$(2,n)$-family $F$ with $N(F)=M(t,n)$, and he conjectured that
$M(t,n)=U(t,n)$ and $F^*$ is the only $(t,n)$-family $F$ with
$N(F)=M(t,n)$ for all $t\ge 3$. In 1992, Leung and Wei
\cite{leung1992} claimed that they proved the above conjecture by
means of a comparison theorem for permanents. But Leung and Wei's
proof has a fatal mistake (see \cite{chang1996}). Hence, the
conjecture is still open. In this paper, we  solve the conjecture.
In fact, we get a more general result for so-called valued
$(t,n)$-family. In what follow, we assume that $t\ge 2$.

 For a sequence of positive integers $(a_1, \cdots, a_n)$, a family
$F=(A_1,\cdots, A_n)$ is called a {\it valued $(t,n)$-family with
valuation $(a_1,\cdots,a_n)$} if $|A_i| = a_i + t$ and
$|\bigcup\limits_{i \in I}A_i|\ge \sum\limits_{i \in I}a_i + t$ for
any $|I| \ge 2$. Note that a $(t,n)$-family $F=(A_1,\cdots, A_n)$
with $N(F)=M(t,n)$ must have $|A_i|=t+1$ for $1\le i\le n$ (see
Lemmas $1$ and $2$ in \cite{chang1989}). Hence, a
 $(t,n)$-family $F$ with $N(F)=M(t,n)$ is a valued $(t,n)$-family with
 valuation $(1,\cdots, 1)$. Let $\bar{F}$ be a valued $(t,n)$-family with valuation $(a_1,
\cdots, a_n)$ satisfying $|\bigcap\limits_{i\in I}A_i| = t$ for any
$|I| \ge 2$. Hence, $F^*$ is $\bar{F}$ with valuation $(1, \cdots,
1)$. Define $M^{'}(t,n,a_1,\cdots, a_n) = \min\{N(F)~|~F$ is a
valued $(t, n)$-family with valuation $(a_1,\cdots,a_n)\}$, and let

$$U^{'}(t,n,a_1,\cdots,a_n)=N(\bar{F})=\sum\limits_{j=0}^{t}\left[{t\choose j}j!\sum\limits_{1\le
i_1< \cdots< i_{n-j}\le n}a_{i_1}\cdots a_{i_{n-j}}\right].$$
 In this paper, we will prove that $M^{'}(t, n,a_1,\cdots,a_n) = U^{'}(t, n,a_1,\cdots,a_n)$
 and $\bar{F}$ is the only valued $(t,n)$-family $F$ with valuation $(a_1,\cdots,a_n)$ satisfying
$N(F) = M^{'}(t, n,a_1,\cdots,a_n)$ for  $t\ge 2$. The conjecture of
Chang \cite{chang1989} is a direct corollary of the conclusion.

Some notations are needed. Suppose $F$ is  a valued $(t,n)$-family
with valuation $(a_1, \cdots, a_n)$.  Let $N=\{1,2,\cdots, n\}$ and
$\mathcal {B}=$ $\bigcup\limits_{i\in N}A_i$, and let  $I_x = \{i\in
N~|~ x \in A_i\}$ and $I_x^c=N-I_x$ for $x \in \mathcal {B}$. The
{\it degree of $x$}, denoted by $\deg x$, is $|I_x|$. A pair of
elements $\{x,y\} \subseteq \mathcal {B}$ is {\it exclusive} if
$I_x\cap I_y^c \neq \emptyset$ and $I_y \cap I_x^c \neq \emptyset$.
An exclusive pair $\{x,y\}$ is {\it saturated} if there exists a
subset $I \subseteq N$ satisfying $I\cap I_x\cap I_y = \emptyset$,
$I\cap I_x\cap I_y^c \neq \emptyset$, $I\cap I_x^c\cap I_y \neq
\emptyset$ and $|\bigcup\limits_{i \in I}A_i| = \sum\limits_{i\in
I}a_i + t$; otherwise, we say an exclusive pair $\{x,y\}$ is {\it
unsaturated}.


\section{An exclusive pair $\{x,y\}$ for a valued $(t,n)$-family}

Assume that $F = (A_1, \cdots, A_n)$ is a valued $(t,n)$-family with
valuation $(a_1,\cdots,a_n)$ and a pair of elements $\{x,y\}$ is
exclusive for $F$. Let

\[
 A_i(x,y)=\left\{
\begin{array}{ll}
 A_i-\{x\}\cup \{y\}, & i \in I_x\cap I_y^c; \\
 A_i,                 &\mbox{otherwise.}
\end{array}
\right.
\]

Then we get a new family $F_{y}^{x} = (A_1(x,y), \cdots, A_n(x,y))$,
but   it is possible that  $F_{y}^{x}$ is  not  a valued
$(t,n)$-family with valuation $(a_1,\cdots,a_n)$. For any
$I\subseteq N$, by calculating $|\bigcup\limits_{i \in I}A_i|$ and
$|\bigcup\limits_{i \in I}A_i(x,y)|$, we can get the relationship
between the two values as follows:

\[
 |\bigcup\limits_{i \in I}A_i(x,y)|=\left\{
 \begin{array}{ll}
 |\bigcup\limits_{i \in I}A_i| - 1, & I\cap I_x\cap I_y=\emptyset, I\cap I_x\cap I_y^c \neq\emptyset,
                                             I\cap I_x^c\cap I_y\neq  \emptyset;\\
 |\bigcup\limits_{i \in I}A_i|, &\mbox{otherwise.}
\end{array}
\right.
\]

Hence, $F_y^x$ is also  a valued $(t,n)$-family with valuation
$(a_1,\cdots,a_n)$ if and only if   $\{x,y\}$ is unsaturated  for
$F$. Furthermore, we have

\begin{theorem}\label{th1}
 A valued $(t,n)$-family with valuation $(a_1,\cdots,a_n)$ satisfying
 $N(F)=M'(t,n,a_1,$ $\cdots,a_n)$ does not contain any unsaturated pair $\{x,y\}$.
\end{theorem}
\begin{proof}
Suppose to the contrary that $\{x,y\}$ is unsaturated for $F$. Then,
$F_y^x$ is also  a valued $(t,n)$-family with valuation
$(a_1,\cdots,a_n)$. We will prove that $N(F_y^x) < N(F)$ and hence
leads to a contradiction.

 Without lose of generality, we can assume that $I_x\cap I_y^c=\{1, \cdots, k_1\}
\neq \emptyset$, $I_y\cap I_x^c=\{k_1+1, \cdots, k_2\} \neq
\emptyset$, $I_x \cap I_y=\{k_2+1, \cdots, k_3\}$ and $I_x^c\cap
I_y^c=\{k_3+1, \cdots,  n\}$. So $F_y^x = (A_1(x,y), \cdots,
A_n(x,y)) = (A_1-\{x\}\cup \{y\}, \cdots, A_{k_1}-\{x\}\cup \{y\},
A_{k_1+1}, \cdots, A_n)$. Let $(x_1,\cdots,x_n)$ be an SDR of
$F_y^x$. Define a function $f$ from the set of all SDRs of $F_y^x$
to the set of all SDRs of $F$ as follows:

(a) if $x_i=y$ for some $i\in\{1,\cdots,k_1\}$ and $x_j=x$ for some
$j\in\{k_2+1,\cdots,k_3\}$, then
\[
(x_1, \cdots, y, \cdots, x, \cdots,  x_n)\rightarrow
 (x_1, \cdots, x, \cdots, y, \cdots,  x_n).
\]

(b) if $x_i=y$ for some $i\in\{1,\cdots,k_1\}$ and $x_j\neq x$ for
all $x_j$, then

\[
(x_1, \cdots, y,   \cdots, x_n)\rightarrow (x_1, \cdots, x, \cdots,
x_n).
\]

(c) otherwise,
\[
(x_1,\cdots,x_n)\rightarrow (x_1,\cdots,x_n).
\]

$f$ is clearly one to one. Define
$$F^{'}=(A_2-\{x,y\},\cdots,A_{k_1}-\{x,y\},A_{k_1+2}-\{x,y\},\cdots, A_n-\{x,y\}).$$

When $t\ge 2$, $F^{'}$ satisfies the Hall's condition  and has an
SDR $(x_2,\cdots,x_{k_1},x_{k_1+2},\cdots,x_n)$. Hence,
 $F$ has an SDR  such as
$$(x,x_2,  \cdots, x_{k_1}, y, x_{k_1+2}, \cdots,  x_n),$$
which is not an $f$-image of an SDR of $F_y^x$, so $f$ is not
subjective. Hence, $N(F_y^x) < N(F)$.
\end{proof}


\section{Saturated pairs of a valued $(t,n)-$family}


For the set $N=\{1, \cdots, n\}$, we define a relation $``\sim''$ on
 $N$ as follows: $i\sim j$ if and only if there exists a subset $I$ satisfying
$\{i, j\} \subseteq I \subseteq N$ and $|\bigcup\limits_{s \in
I}A_s| = \sum\limits_{s \in I}a_s + t$. We claim that $``\sim''$ is
an equivalent relation on $N$. It is obvious that $``\sim''$ is
reflexive and symmetric. If $i\sim j$ and $j\sim k$, then there
exist $I$ and $J$ satisfying $\{i, j\} \subseteq I$,
$|\bigcup\limits_{s \in I}A_s| = \sum\limits_{s \in I}a_s + t$ and
$\{j, k\} \subseteq J$, $|\bigcup\limits_{s \in J}A_s| =
\sum\limits_{s \in J}a_s + t$, respectively. Note that  $ I\cap
J\neq \emptyset$ as $j\in I\cap J$. Hence, we have

\begin{eqnarray*}
\sum\limits_{s \in I\cup J}a_s + t & \le & |\bigcup\limits_{s \in
I\cup J}A_s| = |(\bigcup\limits_{s \in I}A_s) \cup (\bigcup\limits_{s \in J}A_s)|\\
& \le & |\bigcup\limits_{s \in I}A_s| + |\bigcup\limits_{s
\in J}A_s| - |\bigcup\limits_{s \in I\cap J}A_s|\\
& \le & \sum\limits_{s \in I}a_s + t + \sum\limits_{s \in J}a_s + t
- (\sum\limits_{s \in I\cap J}a_s + t) \\
& = & \sum\limits_{s\in I\cup J}a_s + t.
\end{eqnarray*}

So we know that $|\bigcup\limits_{s \in I\cup J}A_s| =
\sum\limits_{s \in I\cup J}a_s + t$ and $\{i, k\} \subseteq I\cup
J$. It implies that  $i\sim k$ and $``\sim''$ is transitive. Hence,
$``\sim''$ is an equivalent relation. So we can classify $N$ into
different classes: $C_1, \cdots, C_m$. If an index set $I\subseteq
N$ satisfies $|\bigcup\limits_{i \in I}A_i| = \sum\limits_{i \in
I}a_i + t$, by the definition of $``\sim''$, we know that
$I\subseteq C_i$ for some $i\in\{1,\cdots,m\}$.

\begin{theorem}\label{th2}
For a valued $(t,n)$-family $F$ with valuation $(a_1, \cdots, a_n)$,
denote by $NSP(F)$ the number of saturated pairs of $F$, then
$NSP(F) \le \sum\limits_{1\le i< j\le n}a_ia_j$.
\end{theorem}
\begin{proof}
We use induction on $n$.  When $n=2$, the conclusion is obvious.

If $|\mathcal {B}| > \sum\limits_{i=1}^{n}a_i + t$, then by the
classification of $N$ under the equivalent relation $``\sim''$, we
get several classes $C_1, \cdots, C_m$ and
 $m\ge 2$.  Without lose of generality, we can assume that $C_1=\{1,
\cdots, k_1\}, \cdots, C_m=\{k_{m-1}+1, \ldots, n\}$. We  get $m$
subfamilies $F_1, \cdots, F_m$ with index sets $C_1, \cdots, C_m$,
respectively. According to the preparation before Theorem \ref{th2},
we know that each saturated pair of $F$ must be saturated for some
subfamily $F_i$. Hence, $NSP(F) \le NSP(F_1) + \cdots + NSP(F_m)$.
By induction,
$$
NSP(F) \le \sum\limits_{1\le i< j\le k_1}a_ia_j + \cdots +
\sum\limits_{k_{m-1}+1\le i< j\le n}a_ia_j < \sum\limits_{1\le
i<j\le n}a_ia_j.
$$

Now  we assume that $|\mathcal {B}| =\sum\limits_{i=1}^{n}a_i + t$.
Let $I$ be an index set satisfying the following conditions: (1)
$|I| \ge 2$; (2) $|\bigcup\limits_{i \in I}A_i| = \sum\limits_{i \in
I}a_i + t$; (3) For $J\subset I$, if $|J| \ge 2$, then
$|\bigcup\limits_{i \in J}A_i|
> \sum\limits_{i \in J}a_i + t$. Since
$|\mathcal {B}| = \sum\limits_{i=1}^{n}a_i + t$, the existence of
such $I$ holds. Now we use different methods to discuss  two cases
$I\subset N$ and $I=N$.

For $I\subset N$, without lose of generality, we can assume that
$I=\{k+1, \cdots, n\}, k\ge 1$. Let $B_1=A_1, \ldots, B_k=A_k,
B_{k+1}=\bigcup\limits_{i=k+1}^{n}A_i$, then $G=(B_1, \cdots,
B_{k+1})$ is a valued $(t,k + 1)$-family with valuation $(a_1,
\cdots, a_k, \sum\limits_{i=k+1}^{n}a_i)$. Let $\{x,y\}$ be an
arbitrary saturated pair  for $F$. There are three subcases: (1)
$\{x,y\}$ is saturated for the subfamily $(A_1, \cdots, A_k)$;
(2) $\{x,y\}$ is saturated for the subfamily $(A_{k+1}, \cdots,
A_n)$; (3) $\{x,y\}$ is unsaturated for both $(A_1, \cdots, A_k)$
and $(A_{k+1}, \cdots, A_n)$. It is easy to see that  $\{x,y\}$ in
the subcase (1) is also saturated for the family $G$.

We claim that  $\{x,y\}$ in the subcase (3) is also  saturated for
$G$. Since $\{x,y\}$ is saturated for $F$ and unsaturated for both
$(A_1, \cdots, A_k)$ and $(A_{k+1}, \cdots, A_n)$, there exist
$\emptyset \neq I_1 \subseteq \{1, \cdots, k\}$ and $\emptyset \neq
I_2 \subseteq I = \{k+1, \cdots, n\}$ such that
$|\bigcup\limits_{i\in I_1\cup I_2}A_i| = \sum\limits_{i\in I_1\cup
I_2}a_i + t$ and  $(I_1\cup I_2)\cap I_x\cap I_y= \emptyset$,
$(I_1\cup I_2)\cap I_x\cap I_y^c\neq \emptyset$, $(I_1\cup I_2)\cap
I_y\cap I_x^c\neq \emptyset$.
Since $|\bigcup\limits_{i\in I_1\cup I_2}A_i| =
\sum\limits_{i\in I_1\cup I_2}a_i + t$ and $|\bigcup\limits_{i\in
I}A_i| = \sum\limits_{i=k+1}^{n}a_i + t$, using the same discussion
in the proof of transitivity of $``\sim''$, we can show that
$|(\bigcup\limits_{i\in I_1}B_i)\cup B_{k+1}| =
|(\bigcup\limits_{i\in I_1}A_i)\cup (\bigcup\limits_{i\in I}A_i)|
=|(\bigcup\limits_{i\in I_1\cup I_2}A_i)\cup (\bigcup\limits_{i\in
I}A_i)|= \sum\limits_{i\in I_1}a_i + \sum\limits_{i=k+1}^{n}a_i +
t$. Under these circumstances, if $\{x,y\}$ is not a subset of
$B_{k+1}$, then $\{x,y\}$ is saturated for $G$.

Now we will prove that $\{x,y\}$ is not a subset of  $B_{k+1}$ in
two cases: $|I_2|\ge 2$ and $|I_2| = 1$.

If $|I_2|\ge 2$, we claim that $I_2=I$. Suppose to the contrary that
$I_2\subset I$. According to  $I=\{k+1, \cdots, n\}$, we know that
$|\bigcup\limits_{i\in I}A_i| = \sum\limits_{i=k+1}^{n}a_i + t$ and
$|\bigcup\limits_{i\in I_2}A_i| > \sum\limits_{i\in I_2}a_i + t$. So
$|(\bigcup\limits_{i\in I}A_i) - (\bigcup\limits_{i \in I_2}A_i)| <
\sum\limits_{i=k+1}^{n}a_i - \sum\limits_{i\in I_2}a_i$. Hence,

\begin{eqnarray*}
|\bigcup\limits_{i\in I_1\cup I}A_i| & = & |\bigcup\limits_{i\in
I_1\cup I_2}A_i| + |(\bigcup\limits_{i\in I-I_2}A_i) -
(\bigcup\limits_{i\in I_1\cup I_2}A_i)|\\
& \le & |\bigcup\limits_{i\in I_1\cup I_2}A_i| +
|(\bigcup\limits_{i\in I}A_i) - (\bigcup\limits_{i\in I_2}A_i)|\\
& < & \sum\limits_{i\in I_1\cup I_2}a_i + t +
\sum\limits_{i=k+1}^{n}a_i - \sum\limits_{i\in I_2}a_i \\
&= & \sum\limits_{i\in I_1\cup I}a_i + t.
\end{eqnarray*}
It contradicts with the fact that $F$ is a valued $(t,n)$-family
with valuation $(a_1, \cdots, a_n)$. Hence, $I_2=I$.

Now we know that   $(I_1\cup I)\cap I_x\cap I_y = \emptyset$, and
hence $I\cap I_x\cap I_y = \emptyset$. Since $|\bigcup\limits_{i\in
I}A_i| = \sum\limits_{i\in I}a_i + t$ and  $\{x,y\}$ is unsaturated
for the subfamily $(A_{k+1}, \cdots, A_n)$, we have  either $I\cap
I_x\cap I_{y}^{c}= \emptyset$ or $I\cap I_{x}^{c}\cap I_{y}=
\emptyset$.  Furthermore, we have either $I\cap I_x = \emptyset$ or
$I\cap I_y = \emptyset$. Therefore, $B_{k+1}=\bigcup\limits_{i\in
I}A_i$ contains at most one of $x,y$, so $\{x,y\}$ is not a subset
of $B_{k+1}$.

If $|I_2| = 1$, without lose of generality, we can assume that
$I_2=\{k+1\}$.  Since $(I_1\cup I_2)\cap I_x\cap I_y = \emptyset$,
we know that $k+1\notin I_x\cap I_y$, which implies that $A_{k+1}$
contains at most one of $x,y$. Assume that $y\notin A_{k+1}$.
Suppose to the contrary that $\{x,y\}$ is a subset of  $B_{k+1}$,
then $y\in \bigcup\limits_{i\in I-I_2}A_i$. By the selection of
$I_1$ and $I_2$, we know that $y\in \bigcup\limits_{i\in I_1\cup
I_2}A_i$, and hence $y\notin (\bigcup\limits_{i\in I-I_2}A_i) -
(\bigcup\limits_{i\in I_1\cup I_2}A_i)$. Then,

$$
|(\bigcup\limits_{i\in I-I_2}A_i) - (\bigcup\limits_{i\in I_1\cup
I_2}A_i)| < |(\bigcup\limits_{i\in I}A_i) - A_{k+1}|.
$$

Since $|A_{k+1}| = a_{k+1} + t$ and $|\bigcup\limits_{i\in I}A_i| =
\sum\limits_{i=k+1}^{n}a_i + t$, we know that

$$
|(\bigcup\limits_{i\in I}A_i) - A_{k+1}| =|\bigcup\limits_{i\in
I}A_i|-|A_{k+1}|= \sum\limits_{i=k+2}^{n}a_i.
$$

Therefore,

\begin{eqnarray*}
|(\bigcup\limits_{i\in I_1}A_i)\cup(\bigcup\limits_{i\in I}A_i)| &=&
|(\bigcup\limits_{i\in I_1}A_i)\cup A_{k+1}\cup
(\bigcup\limits_{i\in I-I_2}A_i)|\\
&= & |(\bigcup\limits_{i\in I_1}A_i)\cup A_{k+1}| +
|(\bigcup\limits_{i\in I-I_2}A_i) -
(\bigcup\limits_{i\in I_1\cup I_2}A_i)| \\
&=& \sum\limits_{i\in I_1}a_i + a_{k+1} + t + |(\bigcup\limits_{i\in
I-I_2}A_i) - (\bigcup\limits_{i\in I_1\cup I_2}A_i)|\\
&<& \sum\limits_{i\in I_1}a_i + \sum\limits_{i=k+1}^{n}a_i + t
\end{eqnarray*}

This contradicts with the fact that $F$ is a valued $(t,n)$-family
with valuation $(a_1, \cdots, a_n)$. Hence, $\{x,y\}$ is not a
subset of $B_{k+1}$.

Now we have shown that when $I\subset N$, any saturated pair
$\{x,y\}$ for $F$ is saturated for either $G$ or the subfamily
$(A_{k+1}, \cdots, A_n)$. Therefore,
$$NSP(F) \le NSP(G) + NSP((A_{k+1}, \cdots, A_{n}))$$
by induction, we have
$$
NSP(G) \le \sum\limits_{1\le i < j\le
k}a_ia_j+(\sum\limits_{l=1}^{k}a_l)(\sum\limits_{m=k+1}^{n}a_m)
$$
and
$$ NSP((A_{k+1}, \cdots, A_{n})) \le \sum\limits_{k+1\le i <
j\le n}a_ia_j$$

Hence, $NSP(F) \le \sum\limits_{1\le i< j\le n}a_ia_j$.

When $I=N$, an exclusive pair $\{x,y\}$ is saturated for $F$ if and
only if $I_x\cap I_y = \emptyset$. Let $C = \{\{x,y\} ~|~I_x\cap I_y
= \emptyset\}$.  Then $NSP(F)=|C|$. Now we calculate $|C|$.

For an arbitrary element $z\in \mathcal {B}$, define $C(z) =
\{\{x,z\} ~|~ I_x\cap I_z = \emptyset\}$. It is not difficult to see
that $|C| = \frac{1}{2}\sum\limits_{z\in \mathcal {B}} |C(z)|$ and
$C(z) = \{\{x,z\} ~|~  I_x\cap I_z = \emptyset\} = \{\{x,z\} ~|~
x\notin \bigcup\limits_{i\in I_z}A_i\}$. So,
$$
|C(z)| =|\mathcal {B}| - |\bigcup\limits_{i\in I_z}A_i| \le
\sum\limits_{i\in I_z^c}a_i.
$$

Therefore,
\begin{eqnarray*}
|C| & \le & \frac{\sum\limits_{z\in \mathcal {B}}\sum\limits_{i\in
I_{z}^{c}}a_i}{2} = \frac{\sum\limits_{z\in \mathcal {B}}
(\sum\limits_{i=1}^{n}a_i - \sum\limits_{i\in I_z}a_i)}{2}\\
& = & \frac{(\sum\limits_{i=1}^{n}a_i + t)(\sum\limits_{i=1}^{n}a_i)
- \sum\limits_{ z\in \mathcal {B}}\sum\limits_{i\in I_z}a_i}{2}\\
& = & \frac{(\sum\limits_{i=1}^{n}a_i + t)(\sum\limits_{i=1}^{n}a_i)
-
\sum\limits_{i=1}^{n}(a_i + t)a_i}{2} \\
&=& \sum\limits_{1\le i< j\le n}a_ia_j.
\end{eqnarray*}

\end{proof}

\section{Exclusive pairs of a valued $(t,n)$-family}

\begin{theorem}\label{th3}
 For a valued $(t,n)$-family $F$ with valuation $(a_1, \cdots, a_n)$,
denote by $NEP(F)$ the number of exclusive pairs of $F$, then
$NEP(F)\ge \sum\limits_{1\le i< j\le n}a_ia_j$. $\bar{F}$ is the
only valued $(t,n)$-family $F$ with valuation $(a_1, \cdots, a_n)$
satisfying $NEP(F)= \sum\limits_{1\le i< j\le n}a_ia_j$.
\end{theorem}
\begin{proof}
 We can assume that $n\ge 2$. For an arbitrary element $z\in \mathcal {B}$,
 $\{x,z\}$ is exclusive for $F$ if and only if   $x\in
\bigcup\limits_{i\in I_z^c}A_i$ and $x\notin \bigcap\limits_{i\in
I_z} A_i$. Define $D(z)=\{\{x,z\} ~|~ \{x,z\}$ is exclusive for
$F\}$. Therefore,
$$D(z)=\{\{x,z\} ~|~ x\in \bigcup\limits_{i\in
I_z^c}A_i - \bigcap\limits_{i\in I_z}A_i\}. $$

Let $\mathcal {A}=\{z|\deg z=n \}$ and $D=\{\{x,y\} ~|~ \{x,y\}$ is
exclusive for $F\}$. Note that $D(z)=\emptyset$ if $z\in \mathcal
{A}$. Then,

\begin{eqnarray*}
|D| & =& \frac{1}{2}\sum\limits_{z\in \mathcal
{B}}|D(z)|=\frac{1}{2}\sum\limits_{z\in \mathcal
{B-A}}|D(z)|\\
& = & \frac{1}{2}\sum\limits_{z\in \mathcal
{B-A}}(|\bigcup\limits_{i\in I_z^c}A_i - \bigcap\limits_{i\in
I_z}A_i|).
\end{eqnarray*}

We first assume that  $\deg z\ge 2$ for all $z\in \mathcal {B-A}$.
Then $|I_z|\ge 2$ and hence $|\bigcap\limits_{i\in I_z}A_i|\le t$
for all $z\in \mathcal {B-A}$. Hence,

$$
|D|>\frac{1}{2}\sum\limits_{z\in \mathcal
{B-A}}(|\bigcup\limits_{i\in I_z^c}A_i| - |\bigcap\limits_{i\in
I_z}A_i|)\ge \frac{1}{2}\sum\limits_{z\in \mathcal
{B-A}}\sum\limits_{i\in I_z^c}a_i.  \eqno(*)
$$

We point out that the inequality strictly holds  as $z\in
\bigcap\limits_{i\in I_z}A_i$ and $z\notin \bigcup\limits_{i\in
I_z^c}A_i$.  To calculate $\sum\limits_{z\in \mathcal
{B-A}}\sum\limits_{i\in I_z^c}a_i$, we  construct a weighted
bipartite graph $G$ as follows: $V(G)=V_1\cup V_2$, where
$V_1=\mathcal {B-A}$ and $V_2=\{A_1,\cdots,A_n\}$; For $z\in V_1$,
if $z\notin A_i$, then $zA_i\in E(G)$ and the weight of  $zA_i$,
denoted by $w(zA_i)$, is $a_i$. 
 So,

$$\sum\limits_{z\in \mathcal {B-A}}\sum\limits_{i\in
I_z^c}a_i=\sum\limits_{z\in V_1}\sum\limits_{zA_i\in E(G)}w(zA_i)
=\sum\limits_{A_i\in V_2}\sum\limits_{zA_i\in E(G)}w(zA_i).
\eqno(**)
$$

Let $|\mathcal {A}|=a $. Obviously, $a\le t$. Each set $A_i$
contains
 $a_i+t-a$ elements in $\mathcal {B-A}$ and there are at least
 $\sum\limits_{j=1}^{n}a_j+t-a$ elements in $\mathcal {B-A}$. By the construction of $G$,  we
 know that the vertex $A_i$ is incident to at least $\sum\limits_{j=1}^{n}a_j-a_i$
 edges in $G$ and the weight of each edge incident to $A_i$ is
 $a_i$. Therefore,
$$
\sum\limits_{A_i\in V_2}\sum\limits_{zA_i\in E(G)}w(zA_i)\ge
\sum\limits_{i=1}^{n}a_i(\sum\limits_{j=1}^{n}a_j-a_i)=(\sum\limits_{i=1}^{n}a_i)^2-\sum\limits_{i=1}^{n}a_i^2.
\eqno(***)
$$

By above  inequalities $(*)$, $(**)$ and $(***)$, we know that  $
|D|>\sum\limits_{1\le i< j\le n}a_ia_j $  if $\deg z\ge 2$ for all
$z\in \mathcal {B}.$

Now we assume that there exists an element $x$ such that $\deg x =
1$, without lose of generality, we assume that $I_x = \{n\}$. Let $k
= \sum\limits_{i=1}^{n}a_i$. We use induction on $k$.

When $k = 2$, then $n = 2$ and $a_1=a_2 = 1$, the conclusion is
obvious.  Assume that $k\ge 3$. As the conclusion is obvious  when
$n=2$, we may assume that $n\ge 3$.

If $a_n = 1$, let $F_1 = (A_1, \cdots, A_{n-1})$, by induction
hypothesis, $NEP(F_1)\ge \sum\limits_{1\le i< j\le n - 1}a_ia_j$ and
$NEP(F_1)= \sum\limits_{1\le i< j\le n - 1}a_ia_j$ implies that
$F_1$ is $\bar{F}$ with valuation $(a_1, \cdots, a_{n-1})$. It is
obvious that the exclusive pairs of $F_1$ are also exclusive for
$F$. Since $(\bigcup\limits_{i=1}^{n-1}A_i) - A_n =
(\bigcup\limits_{i=1}^{n}A_i) - A_n$, we know that $
|\bigcup\limits_{i=1}^{n-1}A_i - A_n|\ge \sum\limits_{i=1}^{n-1}a_i.
$ Obviously, each element $y$ in $(\bigcup\limits_{i=1}^{n-1}A_i) -
A_n$ is exclusive with $x$ for $F$ and $\{x, y\}$ is different from
any exclusive pair of $(A_1, \cdots, A_{n-1})$. Therefore,
$$
NEP(F)\ge \sum\limits_{1\le i< j\le n - 1}a_ia_j +
\sum\limits_{k=1}^{n-1}a_k = \sum\limits_{1\le i< j\le n}a_ia_j.
$$

When $NEP(F)=\sum\limits_{1\le i< j\le n}a_ia_j$, it implies that
$A_n\cap (\bigcup\limits_{i=1}^{n-1}A_i)=t$ and   $NEP(F) -
NEP(F_1)=\sum\limits_{k=1}^{n-1}a_k$. This requires that $F$ is
$\bar{F}$ with valuation $(a_1, \cdots, a_n)$.

If $a_n\ge 2$, let $F_2 = (A_1, \cdots, A_{n-1}, A_n - \{x\})$,
which is a $(t, n)$-family with valuation $(a_1, \cdots, a_{n-1},
a_n-1)$, by induction hypothesis, $NEP(F_2)\ge \sum\limits_{1\le i<
j\le n - 1}a_ia_j + \sum\limits_{k=1}^{n-1}a_k(a_{n}-1)$ and
$NEP(F_2)= \sum\limits_{1\le i< j\le n - 1}a_ia_j +
\sum\limits_{k=1}^{n-1}a_k(a_{n}-1)$ implies that $F_2$ is $\bar{F}$
with valuation $(a_1, \cdots, a_{n-1},a_n-1)$. Similarly, the
exclusive pairs of $F_2$ are also exclusive for $F$,
$|\bigcup\limits_{i=1}^{n-1}A_i - A_n|\ge
\sum\limits_{i=1}^{n-1}a_i$, and each element $y$ in
$\bigcup\limits_{i=1}^{n-1}A_i - A_n$ is exclusive with $x$ for $F$
and $\{x, y\}$ is different from any exclusive pair of $F_2$.
Therefore,
$$
NEP(F)\ge \sum\limits_{1\le i< j\le n - 1}a_ia_j +
\sum\limits_{k=1}^{n-1}a_k(a_{n}-1) + \sum\limits_{k=1}^{n-1}a_k =
\sum\limits_{1\le i< j\le n}a_ia_j.
$$

Similarly, when $NEP(F)=\sum\limits_{1\le i< j\le n}a_ia_j$, it
implies that  $F_2$ must be $\bar{F}$ with valuation $(a_1, \cdots,
a_{n-1}, a_n - 1)$, and since $I_x = \{n\}$, it is obvious that $F$
is $\bar{F}$ with valuation $(a_1, \cdots, a_n)$.
\end{proof}

\section{The conclusion about $N(F)$}
By Theorem \ref{th1},   \ref{th2} and   \ref{th3}, we can easily
arrive at the following conclusion:
\begin{theorem}\label{th4}
$M^{'}(t,n,a_1,\cdots,a_n) = U^{'}(t,n,a_1,\cdots,a_n)$ and
$\bar{F}$ is the only valued $(t,n)$-family $F$ with valuation
$(a_1, \cdots, a_n)$ satisfying $N(F) = M^{'}(t,n,a_1,\cdots,a_n)$
for $t\ge 2$.
\end{theorem}

Applying Theorem \ref{th4} to  $(t, n)$-family, we immediately prove
the conjecture of  Chang in \cite{chang1989}.


\end{document}